\documentclass{amsart}
\usepackage[utf8]{inputenc}
\usepackage{amsmath}
\usepackage{amsfonts}
\usepackage{amssymb}
\usepackage{amsthm}
\usepackage{hyperref}
\usepackage{algorithm}
\usepackage[margin=1.5in]{geometry}
\usepackage{algpseudocode}
\usepackage{etoolbox}
\apptocmd{\sloppy}{\hbadness 10000\relax}{}{}
\usepackage[title]{appendix}
\usepackage{tabu}
\allowdisplaybreaks
\usepackage [english]{babel}
\usepackage [autostyle, english = american]{csquotes}
\MakeOuterQuote{"}

\newtheorem{thm}{Theorem}[section]
\newtheorem{definition}[thm]{Definition}
\newtheorem{lemma}[thm]{Lemma}

\numberwithin{equation}{section}

\author{Kevin G. Hare}

\author{Philip W. Hodges}

\thanks{Research of K. G. Hare was supported by NSERC Grant 2014-03154}
\thanks{Research of P. W. Hodges was supported by NSERC Grant 2014-03154 and the President's Research Award, Faculty of Mathematics, University of Waterloo.}

\title{Applications of Integer Semi-Infinite Programming to the Integer Chebyshev Problem}

\begin{document}

\begin{abstract}
We consider the integer Chebyshev problem, that of minimizing the supremum norm over polynomials with integer coefficients on the interval $[0,1]$. We implement algorithms from semi-infinite programming and a branch and bound algorithm to improve on previous methods for finding integer Chebyshev polynomials of degree $n$. Using our new method, we found 16 new integer Chebyshev polynomials of degrees in the range 147 to 244.

\end{abstract}

\maketitle

\section{Introduction}

\subsection{The Integer Chebyshev Problem}
The supremum norm of a polynomial $p$ over an interval $I$ is defined as
\begin{equation}
\|p(x)\|_I=\sup_{x\in I}|p(x)|
\end{equation}
Let $\mathbb{Z}_n[x]$ denote the polynomials of degree at most $n$ with integer coefficients. The integer Chebyshev problem is the problem of finding a polynomial in $\mathbb{Z}_n[x]$ of  minimal supremum norm on the interval $I$, normalized by the degree of the polynomial; most commonly, and in this paper, the case $I=[0,1]$ is considered.

More precisely, we have the following definition

\begin{definition}
For $n>0$, define 
\begin{equation}
t_{\mathbb{Z},n}(I) = \min_{\substack{p\in \mathbb{Z}_n[x] \\ p \neq 0}}\|p(x)\|^{1/n}_I
\end{equation}
Any degree $n$ polynomial $p\in \mathbb{Z}_n[x]$ such that $\|p(x)\|^{1/n}_I=t_{\mathbb{Z},n}(I)$ is called an integer Chebyshev polynomial of degree $n$. Then the value
\begin{equation}
t_{\mathbb{Z}}(I)=\lim_{n\rightarrow \infty} t_{\mathbb{Z},n}(I)
\end{equation}
is called the integer Chebyshev constant, or integer transfinite diameter for the interval $I$.
\end{definition}

We see that the limit in the definition exists, since
\begin{equation} \label{eq:4}
t_{\mathbb{Z},n+m}(I)^{n+m}\leq t_{\mathbb{Z},n}(I)^nt_{\mathbb{Z},m}(I)^m.
\end{equation}
This follows from the fact that if $p(x)$ and $q(x)$ are integer Chebyshev polynomials of degree $n$ and $m$ respectively, then $t_{\mathbb{Z},n+m}(I)\leq \|p(x)q(x)\|^\frac{1}{n+m}_I$.

We note here that for intervals of size greater than 4, the problem as stated here is solved, with the integer Chebyshev constant for the interval as 1, and the polynomial identically 1 on that interval. 

Let $p(x) = a_n x^n + \dots + a_0$ be a polynomial with integer coefficients, and all of whose conjugates are in an interval $I$.
Then it is relatively easy to show, by means of resultants, that $t_{\mathbb Z}(I) \geq \frac{1}{a_n^{1/n}}$.
See for example \cite{b02}.
By considering the numerator of the iterates of the function 
\[ u(x) = \frac{z (1-z)}{1-3 z (1-z)} \]
we get an infinite family of polynomials, known as the Gorshkov-Wirsing polynomials, all of whose conjugates are in $[0,1]$.
This is used to show that $t_{\mathbb Z}([0,1]) \geq \frac{1}{2.3768417062} \approx 0.4207263$.
Here this constant is explicitly computable to an arbitrary number of digits.
See \cite{g59, m94} for details.
In \cite{be95} it was shown that the lower bound coming from this infinite family is in fact not best possible.
That is, there exists an $\epsilon > 0$ such that $t_{\mathbb Z}([0,1]) \geq 0.4207263\dots + \epsilon$.
At the time no non-trivial lower bound for $\epsilon$ was determined.
Pritsker showed in \cite{p05}, by means of weighted potential theory, that 
    $t_{\mathbb Z}([0,1]) \geq 0.4213$.
Generalizations of these Gorshkov-Wirsing polynomails were considered in \cite{h11}.

Given the submultiplicative nature of $t_{\mathbb Z, n}(I)$ we have
    $t_{\mathbb Z}(I) \leq t_{\mathbb Z, n}(I)$ for all $n$.
This gives a simple method to find an upper bound for $t_{\mathbb Z}(I)$; find large degree polynomials with 
    small supremum norm.
In \cite{be95} a set of $9$ polynomials $p_i(x)$ and exponents $a_i$ were found such that the resulting 
    polynomial $P(x) = p_1(x)^{a_1} \dots p_9(x)^{a_9}$ had small supremum norm.
This was used to show that $t_{\mathbb Z}([0,1]) \leq \frac{1}{2.3605\dots} \approx 0.42364$.
This was done computationally by means of a simplex search on a very large grid of points in $[0,1]$.
The upper bound was improved by Habsieger and Salvy in \cite{hs97} to give $t_{\mathbb Z}([0,1]) \leq 0.4234794$.
In addition, Habsieger and Salvy explicitly computed integer Chebyshev polynomials of degree $n$ for all $n \leq 75$.
Wu in \cite{w03} extended this list of integer Chebshev polynomials up to degree 100.
In \cite{p05,p99} Pritsker introduced techniques from weighted potential theory to help improve 
    lower bounds for the multiplicity of factors in an integer Chebyshev polynomial of degree $n$.
With this he improved the upper bound to $t_{\mathbb Z}([0,1]) \leq 0.4232$.
Flammang in \cite{f09}, by means of auxiliary functions constructed with $28$ polynomials, 
    improved this upper bound to 0.42291334.
In the 2009 PhD thesis of Meichsner \cite{m09}, extensive computational work was done on this problem.
An upper bound of $1/2.36482727 \approx 0.42286386$ was found.
In addition, Meichsner extended the collection of known integer Chebyshev polynomials completely up to degree
    $145$, and also for a large collection of polynomials up to degree $230$. 
It is upon this result that we improve.
Later Flammang, in \cite{f14} improved these upper bounds again to 0.422685 
    by finding a good lower bound in the absolute length of a polynomial.
This last bound is the best known bound.
We unfortunately were not able to improve upon this bound in this paper.

There have been many variations of this problem considered.
For example, the integer Chebyshev problem on sets other than $[0,1]$ has been 
    considered in \cite{ap07, frs97, frs06}.
The monic integer Chebyshev polynomial, denoted $t_M(I)$, was introduced in \cite{bpp03}.
What is interesting in the monic case is that often exact value can be computed. 
In \cite{hs06} it was shown, for all $a$, that $t_M([a, a+1]) = 1/2$.
It was shown that $t_M([0,x])$ is continuous in $x$ for $x > 0$ in \cite{h08b}.
The multivariate case has been looked at in \cite{bp08, m09}.
Applications of this to the leading coefficient in Schur's problem were studied in \cite{p12}.

In this paper, we will apply semi-infinite programming techniques to the problem of finding integer Chebyshev polynomials of degree $n$ on $[0,1]$. We prove 16 polynomials on $[0,1]$ to be integer Chebyshev polynomials of degree 147, 149, 152, 153, 154, 158, 175, 191, 194, 198, 202, 236, 238, 239, 241 and 244 respectively. These new polynomials were not previously known to be integer Chebyshev polynomials, and the previous largest-known integer Chebyshev polynomial was of degree 230. There still remain 20 integer Chebyshev polynomials of degree less than 244 that are unknown.

\subsection{Semi-Infinite Programming}
A semi-infinite programming (SIP) problem is an optimization problem that, in primal form, can be formulated as
\begin{equation*}
\begin{array}{llll}
(P)  &\displaystyle\min_{x\in D} f(x), & D = \{x\in \mathbb{R}: g(x,t) \geq 0, &\forall t\in T\}\\
\end{array}
\end{equation*}
where $x\in \mathbb{R}^n$ and $T$ is an infinite set. Thus the optimization problem is over a finite number of variables subject to an infinite number of constraints.  In a linear semi-infinite programming (LSIP) problem, the objective function $f$ and constraint functions $g(x,t)$, $t\in T$ are affine in the variable $x$. For a general reference on linear semi-infinite programming, see \cite{gl98}.

One of the difficulties of solving LSIP problems is that in order to check the feasibility of a solution $\bar{x}$, we must be able to check whether $\bar{x}\in D$, that is, whether the solution to
\begin{equation*}
(Q) \quad \text{min} \quad g(\bar{x},t) \ \forall t\in T
\end{equation*}
is greater than 0. This is known as the lower-level problem. The problem $(Q)$ is, in general, a non-linear global optimization problem, so may be very difficult to solve. Handling this part is often left to the specific application.

The cutting plane algorithm used in this paper to solve the LSIP problems that will arise is in the class of discretization algorithms. The common theme of this class is to replace the set $T$ by a finite subset $T_k \subset T$, so that we can solve the resulting LP, denoted by $P(T_k)$, with finite constraints, using the simplex method. These algorithms are covered in a general setting in chapter 10 of \cite{gl98}.

The cutting plane algorithm is an iterative algorithm, which produces the next discretized set $T_{k+1}$ by using the solution to the discretized problem on iteration $k$. In particular, given a solution $\overline{x}$ to $P(T_k)$, we solve the lower level problem $(Q)$ to determine a set of values $T'={t\in T : g(\overline{x},t) > 0}$. We then set $T_{k+1}=T\cup T'$ and iterate. Termination occurs when $\overline{x}$ is feasible to within some tolerance $\epsilon$. The constraints given by $t\in T'$ are called cutting planes, as they separate the current solution from some set containing all feasible solutions. This algorithm will be further described in the context it is used in Section 2.

\section{Integer Chebyshev Polynomials of Degree $n$}

The best known methods for finding an integer Chebyshev polynomial $p_n$ for given $n$ use 3 steps:
\begin{enumerate}
\item Determine an initial upper bound for $t_{\mathbb{Z},n}[0,1]$. The typical way to do this is to make use of Equation \eqref{eq:4}, and take
\begin{equation}
c_n=\displaystyle\min_{k\in \{1..n-1\}} \|p_{k}(x)p_{n-k}(x)\|_{[0,1]}
\end{equation}
where $p_i$ is an integer Chebyshev polynomial for each $i\in \{1,..,n-1\}$.
\item Using a set of known factors of integer Chebyshev polynomials of degree $k < n$, determine which are necessary factors of an integer Chebyshev polynomial of degree $n$ with supremum norm less than $c_n$. At the end of this step, we refer to the product of the known factors as $F(x)$ and the product of the unknown factors as $G(x)$.
\item Compute $G(x)$.
\end{enumerate}
Step 3 begins where the techniques of Step 2 fail to produce additional factors, and consists of performing some type of exhaustive search of the possible factors $G(x)$ and choosing the one that minimizes $\|G(x)F(x)\|_{[0,1]}$. Thus $p_n(x)=G(x)F(x)$ is an integer Chebyshev polynomial of degree $n$.

We note the following lemmas, whose proofs can be found in \cite{hs97}:

\begin{lemma}
If $p_n(x)$ is an integer Chebyshev polynomial of degree $n$ for the interval $[0,1/4]$, then $p_n(x(1-x))$ is an integer Chebyshev polynomial of degree $2n$ for the interval $[0,1]$.
\end{lemma}

It follows that there is a symmetric integer Chebyshev polynomial of degree $n$ on $[0,1]$ for all even $n$. We also have

\begin{lemma}
If $n$ is odd, then there exists a symmetric integer Chebyshev polynomial of degree $n$ of the form $(2x-1)p(x(1-x))$
\end{lemma}

These results allow us to work on the interval $[0,1/4]$ instead of $[0,1]$. This proves useful, for example in Step 3 above, by reducing the degree of the unknown factor $G$ before the exhaustive search.

In attempting to find an integer Chebyshev polynomial $p_n(x)$ of degree $n$, we first use methods to determine necessary factors (and their multiplicities) of $p_n(x)$, given that $\|p_n(x)\|\leq c_n$. For all current methods, once we are no longer able produce new factors, we are left with known factors $F$ where $\deg(F) < \deg(p_n)$, if $n$ is sufficiently large.

Although the focus of this paper is on improvements to the search for the missing factor $G(x)$, we will briefly describe some of the basic methods that allow for determining the necessary factors $F$. In particular, a simple method used in \cite{hs97} and others is as follows. Given an upper bound $c_n$ on $t_{\mathbb{Z}}[0,1]$, we have
$$|p_n(x)|=|G(x)|\cdot |F(x)|\leq c_n$$
for all $x \in [0,1]$. We then may be able to prove that a factor of the form $ax-b$ $(a,b\in \mathbb{Z})$ must divide $G(x)$. If
$$c_n \leq \frac{|F(b/a)|}{a^g},$$
where $g=\deg(G)$, then $a^g|G(b/a)|<1$. But $a^gG(b/a)$ is the resultant of $G(x)$ and $ax-b$, and since both have integer coefficients, the resultant is also integer. Thus we conclude that $a^g|G(b/a)|=0$ and so $ax-b$ divides $G(x)$. We then update $F(x)$ to $F(x)\cdot (ax-b)$ and $G(x)$ to $G(x)/(ax-b)$, so $G(x)$ is still the remaining unknown factor.
Similar techniques can be used for algebraic numbers whose conjugates all lie
    within the interval.

This technique can easily be extended to multiple factors by using Markov's bound on $m$th derivatives, as done in \cite{hs97}. More sophisticated techniques have been developed, ranging from Lagrange interpolation in \cite{hs97} to generalized Muntz-Legendre polynomials in \cite{w03} to applications of the simplex method in \cite{m09}. An overview of these methods as well as the most recent and most effective techniques can be found in \cite{m09}.

For the remainder of this section, we use the results of Lemmas 1 and 2 to work instead on the interval $[0,1/4]$ with the transformation $x\mapsto \frac{1-\sqrt{1-4x}}{2}$. We will use $G(x)$ and $F(x)$ to refer to the unknown and known factors respectively on the interval $[0,1/4]$. This transformation halves the degree of $G(x)$, resulting in nearly half as many coefficients to determine.

The algorithm presented here requires the known factors $F(x)$ as input, since we are only concerned with computing $G(x)$. The current best methods for finding $F$ can be found in \cite{m09}, where an implementation in Maple is provided. All results of this section are based on the known factors $F$ provided by that implementation.

\subsection{Formulation as an Integer Semi-Infinite Program}
We can formulate our problem as the following semi-infinite program (SIP) with integer variables $a_i$ and continuous variable $c$:
\begin{equation*}
\begin{array}{llll}
(IP)&\text{minimize} & c &\\
& \text{subject to}& -c \leq \left|F(x)\right| \displaystyle\sum_{i=0}^g a_ix^i \leq c & \forall x\in [0,1/4] \\
& & a_i \in \mathbb{Z} & \forall i\in \{0,\dots,g\}
\end{array}
\end{equation*}
where $\sum_{i=0}^g a_ix^i=G(x)$, $g=\deg(G)$. We note that the constraints are linear in the coefficients $a_i$, so this is a linear mixed-integer SIP. To handle the integrality constraints, we propose a branch and bound algorithm similar to those used for solving integer linear programs.

\subsection{Branch and Bound}
Branch and bound algorithms for minimizing an arbitrary objective function work by producing a tree of search nodes and maintaining a record of the best known solution; in our case, this will be $G^*(x)$, where $\|G^*(x)F(x)\|=c^*$

The root node of the search tree represents the optimization problem over the full solution space, while each branch involves the partition of the search space into two or more components. By finding a lower bound on the objective value of each node, we can eliminate those nodes whose lower bound is greater than $c^*$. We call this "cutting" or "pruning" the branch.

Branch and bound methods to solve integer programs find lower bounds by solving a relaxed linear program, and branch by adding constraints to the integer variables to produce new nodes partitioning the search space. Our approach here is similar, but the relaxed problem is instead a linear SIP.

Each node $N$ is defined by the set of constraints $C$ that have been added to the variables during branching. In order to get a lower bound on the solutions obtainable from $N$, we solve the SIP $(R)$ given by the relaxation of the integrality constraints of $(IP)$ together with the constraints $C$, to get a solution $(\bar{a}, \bar{c})$ with $\bar{c}\leq \nu(IP)$, where $\nu(IP)$ is the optimal solution to $(IP)$. We can find candidates for new best solutions at very low cost by rounding $\bar{a}$ to get $a^*$. If
\begin{equation}
\left\|\left( \sum_{i=0}^g a_i^* x^i \right)F(x)\right\|\leq c^*
\end{equation}
then we replace $G^*(x)$ by $\sum_{i=0}^g a_i^* x^i$. At the end of the algorithm, we output $G^*(x)$, where $G^*(x)F(x)$ is the integer Chebyshev polynomial of degree $n$.

In order to fully specify the branch and bound algorithm, we must describe three more procedures:
\begin{enumerate}
\item A lower bounding method for each node in the branch and bound tree.
\item A branching protocol.
\item A policy for deciding which node to process next.
\end{enumerate}

\subsubsection{Lower bounding at each node}

To find a lower bound at node $N$, we use a cutting plane algorithm as described in the introduction to solve the LSIP. We use an iterative process beginning with a finite subset $T_0 \subset [0,1]$, and compute a solution $(\bar{a}, \bar{c})$ to the LP given by only considering constraints defined by $x\in T_0$ in the SIP $R$ of $N$. We then let
$$f(x)=\left(\sum_{i=0}^g a_i^* x^i\right)F(x),$$
and set
$$T_1=T_0 \cup \{x\in [0,1] : f(x) > \bar{c}, f'(x)=0 \}.$$
We iterate in this fashion until no such values of $x$ can be found or the changes to $\bar{c}$ between successive iterations is sufficiently small. A lower bound for the node $N$ is then provided by the value of $\bar{c}$ from the last iteration.

\subsubsection{Branching protocol}

Before detailing the branching technique, we first note that we need only consider polynomials $G(x)$ such that $a_g \geq 1$, since we are searching for a degree $g$ polynomial and negation does not change the supremum norm. Experimentally, combining this constraint with constraints on the variable $a_0$ (e.g. $a_0=1$) produced the largest likelihood of branch cutting occurring near the start of the algorithm. This motivates our branching protocol: given a set $C$ of constraints for the node $N$, we branch on the coefficient $a_i$ with the smallest index $i$ where $C$ does not already contain an equality constraint for $a_i$. If there are no constraints for $a_i$ in $C$, then we produce four new nodes $N_1, N_2, N_3, N_4$ with constraint sets
\begin{itemize}
\item $C_1 = C \cup \{a_i = \lceil \overline{a}_i \rceil \}$
\item $C_2 = C \cup \{a_i = \lfloor \overline{a}_i \rfloor \}$
\item $C_3 = C \cup \{a_i \geq \lceil \overline{a}_i \rceil +1 \}$
\item $C_4 = C \cup \{a_i \leq \lfloor \overline{a}_i \rfloor -1 \}$
\end{itemize}
respectively.

If there are already inequality constraints on $a_i$ in $C$, then we produce nodes with constraint sets similarly to the above, only where the new constraints do not conflict with the existing constraints in $C$, always producing at least 2 equality constraints.

This branching is a trade off between branching with equality on all possible values of $a_i$, and binary branching using only inequalities. This method avoids computing bounds on the value of $a_i$ while also reducing the total number of nodes created compared to binary branching.

\subsubsection{Selecting the next node to process}

The selection of the next node to process is determined by the data structure in which the nodes are placed. We used a priority queue whereby nodes with smaller $\overline{c}$ are processed first, yielding a "best first search". The heuristic justification of this choice is that since we are searching for nodes with smaller lower bounds first, we will not spend time on a branch which has a better chance of being cut.

\subsubsection{The algorithm}

Algorithm \ref{alg:bnb} is a pseudo-code description of the branch and bound algorithm, given procedures for finding lower bounds and determining new nodes by branching, as described above

\begin{algorithm}
\caption{Branch and Bound}
\begin{algorithmic}[1]
\Procedure{BranchAndBound}{$n,F,c_0$}
\State $UpperBound \gets c_0$
\State $BestSol \gets \text{empty}$
\State $PriQueue \gets (N_0, 0)$ \Comment{$N_0$ has no variable constraints}
\While{$PriQueue \ \text{not empty}$}
\State $(N,\overline{a}) \gets \text{GetNode}(PriQueue)$
\State $a^* \gets \text{Round}(\overline{a})$
\State $c^* \gets \left\|\left( \sum a_i^* x^i \right)F(x)\right\|$
\If {$c^* \leq UpperBound$}
	\State $UpperBound \gets c^*$
	\State $BestSol \gets a^*$
\EndIf
\State $(N_1,N_2,N_3,N_4) \gets \text{Branch}(N)$
\For{$i$ from 1 to 4}
\State $(\overline{a},\overline{c}) \gets \text{LowerBound}(N_i)$
\If{$\overline{c} < UpperBound$}
\State $\text{AddNode}(PriQueue, \text{node}=(N_i, \overline{a}), \text{priority}=\overline{c})$
\EndIf
\EndFor
\EndWhile
\EndProcedure
\end{algorithmic}
\label{alg:bnb}
\end{algorithm}

We note that this algorithm is similar in spirit to the technique used in \cite{hs97} for computing $G(x)$, although described in a different framework. Improvements in our algorithm come from eliminating the need to find bounds on each coefficient and that the branch and bound algorithm is more likely to find good solutions sooner by use of the priority queue.

\subsection{Resultant Search}
We describe here the method used in \cite{m09} to find $G(x)$. This method produced the best results before this paper, and will be used in combination with our branch and bound algorithm to find new integer Chebyshev polynomials.

The main idea is to write the coefficients of $G(x)$ in terms of the resultants of $G(x)$ with $g+1$ linear polynomials, and to perform a search on these resultants, which are also integers, instead of the coefficients themselves. Then we can determine a set of congruences the resultants must satisfy, greatly reducing the number of possible combinations of coefficients.

Consider the $g+1$ linear polynomials $v_ix-w_i$. Then the resultants $r_i=v_i^gG(w_i/v_i)$ for each $i$ give a system of equations of the resultants in terms of the coefficients $a_k$ of $G(x)$. We can solve this system for the variables $a_k$ to get
$$\frac{1}{m_k}\sum_{i=1}^{g+1}t_{k,i}r_i=a_k$$
with $m_k,t_{k,i}\in \mathbb{Z}$. Then
$$\sum_{i=1}^{g+1}t_{k,i}r_i \equiv 0 \ (\text{mod}\  m_k).$$
Letting $M=\text{lcm}(m_1, m_2, ...,m_{g+1})$ gives
$$\sum_{i=1}^{g+1}\frac{M}{m_k}t_{k,i}r_i \equiv 0 \ (\text{mod} \ M).$$
We write this in matrix form as
$$S\mathbf{r}=
\begin{bmatrix}
    s_{1,1} & s_{1,2} & \dots  & s_{1,g+1} \\
    s_{2,1} & s_{2,2} & \dots  & s_{2,g+1} \\
    \vdots & \vdots & \ddots & \vdots \\
    s_{g+1,1} & s_{g+1,2} & \dots  & s_{g+1,g+1}
\end{bmatrix}
\begin{bmatrix}
r_1 \\ r_2 \\ \vdots \\ r_{g+1}
\end{bmatrix}
\equiv
\begin{bmatrix}
0 \\ 0 \\ \vdots \\ 0
\end{bmatrix}
\ (\text{mod} \ M).
$$

Using a form of Gaussian elimination with the Euclidean algorithm, we can reduce the matrix $S$ to an upper triangular matrix $S'$. Given bounds on each $r_i$, we can then try all possible values of $\mathbf{r}$ by back substitution in $S'$, solving the linear congruence in one variable at each step. We find bounds on the $r_i$ by noting that $r_i=v_i^g|G(w_i/v_i)|\leq v_i^gc_n/|F(w_i/v_i)|$ or by solving the SIP
$$
\begin{array}{llll}
    &\text{minimize}  & r_i=v_i^g|G(w_i/v_i)| &\\
	&\text{subject to} &\left| G(x) \right| \left|F(x)\right| \leq c_n &\forall x\in [0,1/4]\\
\end{array}
$$
with the coefficients of $G(x)$ as variables, using some subset of $[0,1/4]$ for the constraints (around 200 constraints), whichever gives a better bound.

For each valid $\mathbf{r}$, we can compute $G(x)$ and find the $\mathbf{r}$ giving $G(x)$ that minimizes $\|G(x)F(x)\|$. Then $G(x)F(x)$ is an integer Chebyshev polynomial of degree $n$ on $[0,1/4]$, which we can then transform to find the symmetric integer Chebyshev polynomial on $[0,1]$.

\subsection{Combining Methods}
Along with the other methods given in \cite{m09}, the resultant search method allowed the computation of all integer Chebyshev polynomials up to degree 145, and many up to degree 230, and included all integer Chebyshev polynomials for which the degree of the missing factor $G(x)$ on the interval $[0,1/4]$ was less than 15. Unfortunately, the exponential nature of the branch and bound method prevents scaling to the same extent. For example, in step 3 of finding the integer Chebyshev polynomial of degree 120, $G(x)$ has degree 12. The resultant search method completes in 3 minutes, while branch and bound takes over 8 hours. However, fast progress by the branch and bound algorithm in the early stages motivates a combination of the two methods wherein we start off branching before doing a resultant search on the remaining variables.

More precisely, we observed that large amounts of branch cutting occurred in the early stages of the branch and bound algorithm. This was especially true when the signs of the early variables did not match with the signs of the variables in the minimizing factor. For example, if $a_0=1$ for the factor $G^*(x)$ minimizing $\|G(x)F(x)\|$, then the node with constraint $a_0=-1$ would always produce a lower bound greater than $c_n$. Such asymmetry does not exist when searching the resultants.

We exploit this early effectiveness of branching by combining the two methods; we first start with branch and bound, and branch until a specified variable is reached. Given the determined coefficients, we then do a resultant search to find the remaining coefficients that minimize $\|G(x)F(x)\|$ for each node.

Suppose that we branch until coefficient $j$. Let $\overline{a}_0,\overline{a}_1,...,\overline{a}_j$ be the known coefficients and let $a_{j+1},a_{j+2},...,a_g$ be the unknown coefficients. Then the resultant of $G(x)$ with $v_ix-w_i$ is
$$r_i=v_i^g\left( \sum_{k=0}^j\overline{a}_k(w_i/v_i)^k+ \sum_{k=j+1}^ga_k(w_i/v_i)^k \right).$$
Rearranging, we have
$$\overline{r}_i=r_i-v_i^g\sum_{k=0}^j\overline{a}_k(w_i/v_i)^k=v_i^g\sum_{k=j+1}^ga_k(w_i/v_i)^k$$
with $\overline{r}_i\in \mathbb{Z}$. We can use the system of equations between the $\overline{r_i}$ and the $a_k$ given in this way for the resultant search, shifting the bounds on $r_i$ by $\sum_{k=0}^j\overline{a}_k(w_i/v_i)^k$ to get bounds on $\overline{r_i}$.

In practice, branching until 11 variables remained proved effective. Branching later resulted in set up computations for the resultant search being wasted on a small search, and branching earlier often made the resultant search too large to be feasible in our implementation.

A running time comparison on a select few examples between the combined method and the branch and bound and resultant search alone is given in Table \ref{table:time}. As can be seen, the combined method is not unilaterally better than the resultant search alone, particularly for smaller $n$. However, for larger $n$, we can see that there are instances where the combined method will finish quite quickly, while the resultant search alone does not finish within any reasonable timeframe.

\begin{table}
\begin{center}
\begin{tabular}{ l | c | c | c }
	$n$ & Branch and Bound & Resultant Search & Combined Method \\
	\hline
	120 & 8.1 hours & 182 seconds & 1560 seconds \\
	145 & - & 17.5 hours & 24.0 hours\\
	154 & - & - & 15.7 hours \\
	199 & - & 87.1 hours & 39.9 hours \\
  	\hline 
\end{tabular}
\caption{Time for algorithms to compute integer Chebyshev polynomial of degree $n$. An entry of "-" indicates the program did not terminate within 7 days. Tests were run on AMD Opteron 6174 2.2 GHz processors.}
\label{table:time}
\end{center}
\end{table}

We note two additional advantages of the combined method that are inherited from the branch and bound algorithm. The "best-first" search of the branch and bound method applies here as well, and the variable sets chosen for the first few searches will be the variable sets indicated to be most promising by the LP solutions at all active nodes. Also, as the vast majority of the computations happen within nodes and not between them, the branching portion of the algorithm is very simple to parallelize. Furthermore, it parallelizes nearly perfectly once the number of active nodes exceeds the number of cores available. Future research could implement a parallelized version of the algorithm that could help find integer Chebyshev polynomials that have proven elusive thus far.

\subsection{Results}
Our implementation in Maple of the combined branch and bound and resultant search method, as well as additional data files, can be found at \cite{kghareweb}. 

Using the combined branch and bound and resultant search algorithm, we were able to find 16 more integer Chebyshev polynomials not given in \cite{m09}, listed in Table \ref{table:ICPs}. The polynomials $h_i$ are given in Appendix A. The missing factors found for these all had degrees in the range 14-16, except for the degree 202 integer Chebyshev polynomial, which had a degree 17 missing factor.

\begin{table}
\begin{center}
\begin{tabular}{ l | c | r }
	$n$ & $p_n $ & $t_{\mathbb{Z}, n}[0,1]$ \\
	\hline			
  	147 & $h_1 ^{48}h_2 ^{17}h_3 ^6h_5 ^2h_{10}h_{14}$ & 0.42591455 \\
	149 & $h_1 ^{47}h_2 ^{17}h_3 ^6h_5 ^3h_{10} h_{14} $ & 0.42578804\\
	152 & $h_1 ^{47}h_2 ^{16}h_3 ^6h_5 ^2h_{10} h_{23} $ & 0.42577465\\
	153 & $h_1 ^{48}h_2 ^{19}h_3 ^5h_4 h_5 ^2h_7 h_{10} h_{14} $ & 0.42547485\\
	154 & $h_1 ^{49}h_2 ^{18}h_3 ^6h_5 ^3h_{10} h_{14} $ & 0.42548736\\
	158 & $h_1 ^{51}h_2 ^{18}h_3 ^6h_5 ^3h_{10} h_{14} $ & 0.42536299\\
	175 & $h_1 ^{56}h_2 ^{23}h_3 ^6h_4 h_5 ^2h_7 h_{10} h_{14} $ & 0.42542222\\
	191 & $h_1 ^{60}h_2 ^{21}h_3 ^8h_5 ^3h_{10} h_{14} h_{15} $ & 0.42512849\\
	194 & $h_1 ^{61}h_2 ^{22}h_3 ^8h_5 ^3h_{10} h_{14} h_{15} $ & 0.42517829\\
	198 & $h_1 ^{63}h_2 ^{22}h_3 ^8h_5 ^3h_{10} h_{14} h_{15} $ & 0.42505003\\
	202 & $h_1 ^{64}h_2 ^{24}h_3 ^9h_4 h_5 ^3h_6 h_{10} h_{14} $ & 0.42514131\\
	236 & $h_1 ^{75}h_2 ^{28}h_3 ^9h_4 h_5 ^4h_{10} h_{14} h_{15} $ & 0.42434377\\
	238 & $h_1 ^{76}h_2 ^{28}h_3 ^9h_4 h_5 ^4h_{10} h_{14} h_{15} $ & 0.42468031\\
	239 & $h_1 ^{76}h_2 ^{27}h_3 ^{10}h_5 ^3h_{10} h_{14} h_{21} $ & 0.42461390\\
	241 & $h_1 ^{77}h_2 ^{27}h_3 ^{10}h_5 ^3h_{10} h_{14} h_{21} $ & 0.42448242\\
	244 & $h_1 ^{78}h_2 ^{30}h_3 ^{11}h_4 h_5 ^4h_6 h_{10} h_{14}$ & 0.42456112\\
  	\hline  
\end{tabular}
\caption{New integer Chebyshev polynomials for the interval $[0,1]$. The polynomials $h_i$ are given in Appendix A.}
\label{table:ICPs}
\end{center}
\end{table}

\section{Conclusion}

In this paper, we presented an improved method for finding integer Chebyshev polynomials of degree $n$. We implemented additions to the standard semi-infinite programming methods, including a branch and bound technique combined with previous methods. These efforts yielded 16 new integer Chebyshev polynomials of degree ranging from 147 to 244. There remain 20 integer Chebyshev polynomials of degree less than 244 that are not known.

In the search for integer Chebyshev polynomials, any improvement on either computing the necessary degree of known factors or on the exhaustive search to find remaining factors will likely yield improved results. As noted previously, the branch and bound algorithm naturally parallelizes, and this would be a good way to continue to find integer Chebyshev polynomials, although there are limitations compared to an algorithmic improvement.

\appendix

\section{Polynomials}

We list here the polynomials $h_i(x)$ used in Section 2. These polynomials form a previously known comprehensive list of all known factors of integer Chebyshev polynomials of degree $n$; that is, we found no new integer Chebyshev polynomials with previously unseen factors.

\begin{flalign*}
h_{1}(x)=\ & x(1-x) \\
h_{2}(x)=\ & 2x-1 \\ 
h_{3}(x)=\ &    5 x^{2}-5 x+1 \\ 
h_{4}(x)=\ &   6 x^{2}-6 x+1 \\ 
h_{5}(x)=\ &   29 x^{4}-58 x^{3}+40 x^{2}-11 x+1 \\ 
h_{6}(x)=\ &   33 x^{4}-66 x^{3}+45 x^{2}-12 x+1 \\ 
h_{7}(x)=\ &   34 x^{4}-68 x^{3}+46 x^{2}-12 x+1 \\ 
h_{8}(x)=\ &   49 x^{4}-98 x^{3}+69 x^{2}-20 x+2 \\ 
h_{9}(x)=\ &   161 x^{6}-483 x^{5}+575 x^{4}-345 x^{3}+109 x^{2}-17 x+1 \\ 
h_{10}(x)=\ &   169 x^{6}-507 x^{5}+601 x^{4}-357 x^{3}+111 x^{2}-17 x+1 \\ 
h_{11}(x)=\ &   181 x^{6}-543 x^{5}+644 x^{4}-383 x^{3}+119 x^{2}-18 x+1 \\ 
h_{12}(x)=\ &   193 x^{6}-579 x^{5}+683 x^{4}-401 x^{3}+122 x^{2}-18 x+1 \\ 
h_{13}(x)=\ &   821 x^{8}-3284 x^{7}+5555 x^{6}-5171 x^{5}+2886 x^{4}-985 x^{3}+200 x^{2}-22 x+1 \\ 
h_{14}(x)=\ &   941 x^{8}-3764 x^{7}+6349 x^{6}-5873 x^{5}+3243 x^{4}-1089 x^{3}+216 x^{2}-23 x+1 \\ 
h_{15}(x)=\ &   961 x^{8}-3844 x^{7}+6478 x^{6}-5980 x^{5}+3291 x^{4}-1100 x^{3}+217 x^{2}-23 x+1 \\ 
h_{16}(x)=\ &   4921 x^{10}-24605 x^{9}+53804 x^{8}-67586 x^{7}+53866 x^{6}-28388 x^{5}+9995 x^{4}\\&-2317 x^{3}+338 x^{2}-28 x+1 \\ 
h_{17}(x)=\ &   31169 x^{12}-187014 x^{11}+502099 x^{10}-796200 x^{9}+828936 x^{8}-595698 x^{7}\\&+302334 x^{6}-108945 x^{5}+27600 x^{4}-4783 x^{3}+537 x^{2}-35 x+1 \\ 
h_{18}(x)=\ &   43609 x^{12}-261654 x^{11}+704777 x^{10}-1125390 x^{9}+1184854 x^{8}-865270 x^{7}\\&+448776 x^{6}-166327 x^{5}+43659 x^{4}-7905 x^{3}+936 x^{2}-65 x+2 \\ 
h_{19}(x)=\ &   161429 x^{14}-1130003 x^{13}+3599830 x^{12}-6908941 x^{11}+8913112 x^{10}-\\&8165339 x^{9}+5470288 x^{8}-2718775 x^{7}+1005970 x^{6}-275399 x^{5}+54846 x^{4}\\&-7697 x^{3}+719 x^{2}-40 x+1 \\ 
h_{20}(x)=\ &   161813 x^{14}-1132691 x^{13}+3608246 x^{12}-6924493 x^{11}+8931952 x^{10}\\&-8181043 x^{9}+5479474 x^{8}-2722543 x^{7}+1007031 x^{6}-275594 x^{5}+54867 x^{4}\\&-7698 x^{3}+719 x^{2}-40 x+1 \\ 
h_{21}(x)=\ &   161929 x^{14}-1133503 x^{13}+3610755 x^{12}-6928991 x^{11}+8937122 x^{10}\\&-8185014 x^{9}+5481532 x^{8}-2723251 x^{7}+1007185 x^{6}-275613 x^{5}+54868 x^{4}\\&-7698 x^{3}+719 x^{2}-40 x+1 \\ 
h_{22}(x)=\ &   887981 x^{16}-7103848 x^{15}+26189139 x^{14}-59006633 x^{13}+90856296 x^{12}\\&-101276631 x^{11}+84454852 x^{10}-53688009 x^{9}+26265936 x^{8}-9911593 x^{7}\\&+2872148 x^{6}-631701 x^{5}+103263 x^{4}-12115 x^{3}+961 x^{2}-46 x+1 \\ 
h_{23}(x)=\ &   907201 x^{16}-7257608 x^{15}+26750188 x^{14}-60243176 x^{13}+92693614 x^{12}\\&-103221560 x^{11}+85965780 x^{10}-54562008 x^{9}+26643715 x^{8}-10032840 x^{7}\\&+2900545 x^{6}-636399 x^{5}+103781 x^{4}-12149 x^{3}+962 x^{2}-46 x+1
\end{flalign*}

\bibliographystyle{plain}
\bibliography{references}

\end{document}